\documentclass[11pt,fleqn,a4paper]{article}
\usepackage{amsmath}
\usepackage{amssymb}
\usepackage{makeidx}
\usepackage{pdfpages}
\usepackage{pgf,tikz}
\usepackage{mathrsfs}
\usetikzlibrary{arrows}
\usepackage{amssymb,amsmath,graphicx}
\usepackage{epstopdf}
\usepackage{multirow}
\usepackage{caption}
\usepackage{blkarray}
\usepackage{setspace}
\usepackage{booktabs}
\usepackage{bigstrut}
\usepackage{needspace}
\usepackage[normalem]{ulem}

\usepackage{titling}
\begin{document}
\sloppy

\newtheorem{axiom}{Axiom}[section]
\newtheorem{example}[axiom]{Example}
\newtheorem{proposition}[axiom]{Proposition}

\newcommand{\eps}{\varepsilon}
\newcommand{\real}{\mathbb{R}}
\newcommand{\nat}{{\mathbb{N}}}
\newcommand{\seq}[1]{\langle #1\rangle}

\newcommand{\crcl}{\mbox{$\mathcal{R}$}}
\newcommand{\bound}{\mbox{$\mathcal{L(R)}$} }
\newcommand{\ccc}{{\cal C}}
\newcommand{\VL}{\mbox{$\mathrm{VL}$}} 
\newcommand{\VR}{\mbox{$\mathrm{VR}$}}


\title{2-Period Balanced Travelling Salesman Problem: a polynomially solvable case and heuristics}

\author{
  Vladimir Deineko
  \thanks{{\tt V.Deineko@warwick.ac.uk}.
  Warwick Business School, Coventry, CV4 7AL, United Kingdom}
  \thanks{Corresponding author; tel.\ +44 02476524501}
  \and Bettina Klinz
  \thanks{{\tt klinz@math.tugraz.at}.
Institut f\"ur Diskrete Mathematik, TU Graz, Steyrergasse 30, 8010 Graz, Austria}
 \and Mengke Wang
  \thanks{{\tt m.wang.36@warwick.ac.uk}.
 Warwick Business School, Coventry, CV4 7AL, United Kingdom}
}
\date{}
\maketitle


\begin{abstract}
 We consider the ${\cal NP}$-hard 2-period balanced travelling
  salesman problem.  In this problem the salesman  needs to visit a
  set of customers in two time periods. A given subset of the customers
  has to be visited in both periods while the rest of the customers need to be visited only once, in any of the two periods. Moreover, it is required that the number of customers   visited in the first period does not differ from the respective number in the
  second period by more than $p$, where $p$ is a given small constant.
  The objective is to find
  tours for visiting the customers in the two periods with minimal total
  length. We show that in the case when the underlying distance matrix is a
  Kalmanson matrix, the problem can be solved in polynomial time. For
  the general case, we propose two new heuristics. A combination of these 
  two heuristics has shown very favourable
  results in computational experiments on benchmark instances
  from the literature.
  
\medskip\noindent{\em Keywords.}
Combinatorial optimization; Kalmanson matrix;  2-period balanced travelling salesman problem; vehicle routing problem; polynomially solvable case. 
\end{abstract}

\section{Introduction}
\label{sec:introduction}
\nopagebreak

\subsection{Definitions and related work}
In the well-known travelling salesman problem (TSP), an $n\times n$ distance matrix $(c_{ij})$ among
$n$ locations is given: location $1$ is where the salesman starts and ends the journey, and the other $n-1$ locations are for the customers to be visited. The objective of the salesman is to find a shortest route for visiting all customers. The reader is referred  to the books \cite{ ABCC, Gutin, TSP} for the wealth of knowledge on this problem.

A feasible solution to the TSP can be represented as a  \emph{tour}
$\tau=\seq{1=\tau_1,\tau_2,\ldots,\tau_n,1=\tau_{n+1}}$ on the set
$\{1,2,\ldots,n\}$. We will
refer to the items in tours interchangeably as \emph{points} or \emph{nodes}. The salesman starts from node $1$, visits node
$\tau_2$, and so on, and eventually returns to the initial node $1$.  
Given a distance matrix $(c_{ij})$, the length of the tour $\tau$ is calculated as $c(\tau)=\sum_{i=1,\ldots,n}c_{\tau_i,\tau_{i+1}}$. We assume that all distance matrices considered in this paper are \emph{symmetric}.

In the 2-period travelling salesman problem (2-TSP) the salesman has to visit a set of customers in two periods. A given subset of the customers has to be visited twice -- once in each of the two periods; while the rest of the customers have to be visited only once, in any of the two periods. The objective is to find tours for the two periods with the minimal total length. 

We assume that the salesman always starts and ends the journeys in location 1. For the 2-TSP this assumption means that location 1 is visited twice, i.e.\ 1 belongs to both tours in any feasible solution.

The 2-TSP can be viewed as a simplification of the $m$-period TSP (see e.g.\  \cite{BPS,Chao,Paletta1,Paletta2}) and also as a model for some practical applications of the vehicle routing problem (VRP).
In particular, Butler, Williams, and Yarrow \cite{ButW97} considered the 2-TSP as a model for a milk collection problem in Ireland. Doerner et al.\ \cite{DGHL08} have described a similar model for the two-day blood delivery service of the Austrian red cross. Hamzadayi, Topaloglu, and  Kose
\cite{HTK} considered  $m$-period TSP as a model for delivery goods by a soft drinks company in Turkey.

In this paper we focus on the \emph{balanced} 2-TSP, where an additional constraint requests that the number of customers visited in each period differs by no more than $p$, where $p$ is a given small constant. Different variants of balanced TSPs can also be found in \cite{GPTV, KSDS, LP}.

\begin{figure}
\unitlength=1cm
\begin{center}
{\begin{picture}(6.5,5.5)
\includegraphics[scale=1.1]{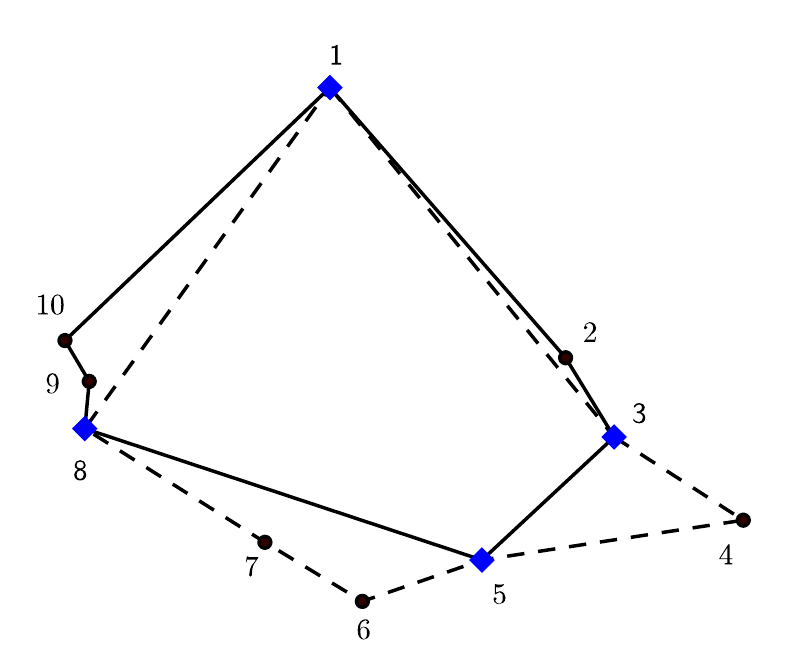}
\end{picture} }
\end{center}
\caption{Illustration to the definition of the balanced 2-TSP; 
locations $1$, $3$, $5$, and $8$ are visited in both periods; the same number of locations is visited in each period. }
\label{fig:2pTSP}
\end{figure}

Bassetto and Mason \cite{BassettoM11} considered a special case of the
balanced 2-TSP where $p=1$ and where the points are located in the Euclidean
plane. 
Fig.~\ref{fig:2pTSP} illustrates a solution to an instance of the  balanced
2-TSP with $10$ points in the Euclidean plane. Four out of ten points in this instance are visited in both periods.  
The set of points shown in Fig.~\ref{fig:2pTSP} happens to lead to a distance
matrix which is a so-called Kalmanson matrix, cf.~\cite{DRVW}.

A \emph{symmetric} $n\times n$ matrix $C$ 
is called a \emph{Kalmanson} matrix if it fulfils the 
\emph{Kalmanson conditions}:
\begin{eqnarray}
c_{ij}+c_{\ell m}\le c_{i\ell}+c_{jm}, \ \ \ \
c_{im}+c_{j\ell}\le c_{i\ell}+c_{jm},\ \ 
\mbox{~for all~} 1\le i<j<\ell<m\le n. \label{kalm1.c}
\end{eqnarray}

Kalmanson~\cite{Kalmanson} noticed that if $n$ points in the Euclidean plane
are located on the boundary of their convex hull and numbered along the convex
hull, then the distance matrix for these points satisfies conditions
(\ref{kalm1.c}). He proved that the TSP with a Kalmanson distance matrix is solved to optimality 
by the tour $\pi=\seq{1,2,3,\ldots,n-2,n-1,n,1}$. Notice that if a distance matrix satisfies conditions
(\ref{kalm1.c}), it does not necessarily mean that the points are on the
boundary of a convex hull (or in the Euclidean plane). For example, the points in the example in
Fig.~\ref{fig:2pTSP}
are not on the boundary of the convex hull, however the distance matrix for
these points is a Kalmanson matrix (for details and the coordinates of the
points see Fig. 1 in \cite{DRVW}).

 It can be easily shown that
cyclic renumbering of rows and columns does not destroy property
(\ref{kalm1.c}): if matrix $C=(c_{ij})$ is a Kalmanson matrix, then the matrix
$C_\sigma=(c_{\sigma(i)\sigma(j)})$ with rows and columns permuted with the
permutation $\sigma=(k,k+1,\ldots,n,1,2,\ldots,k-1)$, is also  a Kalmanson  matrix. Permutation $\sigma$ is called a cyclic shift of the identity permutation. Adding constants (positive or negative) to rows and columns in the distance matrix does not destroy  property (\ref{kalm1.c}) either.

Kalmanson matrices possess nice structural properties that have been used to identify special cases in a wide range of combinatorial optimization problems 
(see \cite{Cela2018,CFSS,DRW,KW,PSW}).
In particular, Kalmanson matrices play a key role in the so-called 
master tour problem \cite{DRW}.

An  \emph{optimal} tour for the TSP is called a \emph{master tour} if after
deleting any subset of points from the tour the so obtained sub-tour is still
optimal for the remaining set of points. Given a distance matrix, \emph{the
  master tour problem} is to find out whether it is possible to construct a
master tour. 

A master tour exists if and only if the underlying distance matrix is a Kalmanson matrix (Theorem 5.1 in \cite{DRW}). In this case the master tour is either the tour $\pi$ or the inverse of it, tour 
$\pi^-=\seq{1,n,n-1,\ldots,2,1}$.

\smallskip

The 2-TSP is related to the more general vehicle routing problem (VRP).
In the VRP a set of given customers are to be visited
by a given fleet of vehicles which are located in one or several depots.
The task is to find an optimal set of routes for the vehicles fulfilling
the requests of the customers while respecting the given constraints
which may vary vastly depending on which VRP variant is considered.
A standard VRP objective is to minimize the total travelled distance.
For further
details the reader is referred to \cite{CAGR,Drexl,Golden,LKH,Toth, VCGP1,VCGP2}.

Note that the distance matrices which result from points located on a line,
a cycle, or a tree  are all Kalmanson matrices. Hence the
papers~\cite{ Archetti2, Hassin, Labbe, YuLiu} can be viewed as
investigations of the VRP for special classes of Kalmanson matrices.

The 2-TSP can be modelled within the VRP framework as a VRP with
only 2 vehicles  
(2-VRP) as follows. One of 
the $m=|S|$ customers that have to be visited in both periods is chosen as a 
depot. For each of the other $m-1$ customers in $S$ an identical copy (i.e. a new
customer in the same location) is created. Every two identical customers are
allocated to two different vehicles. The allocation of these customers to
the vehicles is fixed: each vehicle has to serve $m-1$ fixed customers, while
$k=n-|S|$ customers are to be served by only one of the vehicles.  
In Section~\ref{sec:framework}  we use a VRP related terminology and the 2-VRP model as a useful
tool for representing our algorithms for the 2-TSP. Recall that ``2" in the
abbreviation 2-VRP refers to the number of vehicles and  not to
the capacity of the vehicle as for instance in the abbreviation 2SDVRP for the
Split Delivery Vehicle Routing Problem \cite{Archetti1} (see also \cite{Hassin}).

\subsection{\bf Our contribution and structure of the paper} 
In Section~\ref{sec:solvable2-TSP}, we focus on the balanced 2-TSP with a
Kalmanson distance matrix. We prove that this special case can be solved
in polynomial time by dynamic programming (DP). We then use the DP
recursions in a  heuristic for the 2-TSP on arbitrary
distance matrices.
Computational experiments  on the benchmark problems from the
literature demonstrate that our heuristic which is motivated by a polynomially
solvable case is competitive with previously published algorithms. Thus we demonstrate
that the approach of using efficiently solvable special cases to
devise heuristics for general hard problems could be a fruitful one also
in other settings. This is a major contribution of this paper.

While designing an algorithm for finding an optimal solution for the 2-TSP with a Kalmanson distance matrix, we use a simple methodological approach - ``one sequence - two tours": an optimal solution is found among sequences which combine two tours. In this way the problem of allocating nodes to tours and the problem of optimising the tours are solved at one go. In Section~\ref{sec:framework}, we generalise this approach and adapt the well known Held and Karp \cite{HeldKarp} DP recursions for finding an optimal sequence containing two tours for the 2-TSP with arbitrary distance matrices. The  2-VRP, which is a generalisation of the 2-TSP is then considered. To deal with large size instances of the 2-VRP,  a disassemble-aggregate heuristic is suggested. The heuristic transforms a 2-VRP instance into a collection of small size
instances and solves them by using the DP recursions. Due to the flexibility of the DP technique, various additional constraints can be included in the consideration. 
The new heuristic for the 2-VRP (balanced 2-TSP) is another contribution of this paper.


In Section~\ref{sec:emp} the results of computational experiments are reported. In the experiments, two types of initial solutions have been used: solutions obtained by a heuristic based on the special solvable case and solutions obtained with a random allocation of nodes to  tours. The test results have shown that the approach ``one sequence - two tours" incorporated in the suggested heuristic is very robust to changes in the quality of the initial solutions. The set of 60 benchmark instances introduced in \cite{BassettoM11} have been used in the experiments.  The authors in \cite{BassettoM11} have presented two types of results: (a) results obtained by the algorithms proposed, and (b) results obtained by a
hybrid approach where a computer algorithm was run first and the obtained
results were  then improved by human intervention on the basis of a graphical visualisation of the results. In our experiments, 58 results  of  type (a), and  48 results of type (b) have been improved.  

Section~\ref{sec:sum} provides summary and conclusions.
Technical details and supportive material are placed in the appendices.

\section{Polynomially solvable case of the balanced  2-TSP}
\label{sec:solvable2-TSP}


In this section we consider the \emph{balanced} 2-TSP with an $n\times n$ Kalmanson distance matrix $C$.  Let $S\subset\{1,2,\ldots,n\}$ denote the subset of customers to be visited in both periods.
We refer to these customers as \emph{fixed} nodes.   Given an integer parameter $p$, the objective in the balanced 2-TSP is to find two tours with the minimal total length such that every fixed node 
appears in both tours, every node which is {\em not} in $S$ appears only in one of the tours, and the number of nodes in the tours differs by no more than $p$.  For example, if $p=1$, then for even $n$ the tours in the balanced 2-TSP have to contain the same number of nodes. For arbitrary $p$, each tour has to contain at least
$\lceil(n+|S|-p)/2\rceil$ nodes or equivalently no more
than $\lfloor(n+|S|+p)/2\rfloor$ nodes.
Without loss of generality we assume that $1\in S$ (note that the Kalmanson
property is conserved by a cyclic shift, hence any node can be chosen to be the
start node). Let $s^*$ denote the maximal item in $S$.

The tours $\seq{1,2,\ldots,n-1,n,1}$ and $\seq{1,n,n-1,\ldots,2,1}$ are the
\emph{master} tours for the TSP with a Kalmanson matrix. If a set of indices is removed from the set of rows and the set of columns in a Kalmanson matrix, the so obtained submatrix is still a Kalmanson matrix. 
It means that if an optimal allocation of items to two tours in the 2-TSP is found, an optimal ordering  in     each tour can  be found by sorting the items either in increasing or
decreasing order.  We will search for an optimal solution of the balanced 2-TSP among the solutions 
where the items in the first tour are sorted in increasing order, and in the
second tour in decreasing order. 

Two {\em feasible} tours can be combined into a single feasible sequence $Q$ as follows.
We start with a partially created
sequence $Q=\seq{1,\ldots,1}$ where only the placement
of the start and end node is known. 
There are the following options for placing node 2 in the sequence.
\begin{itemize}
  \item[(a)] If
node 2 is a fixed node, then we place one copy of $2$ after the start node $1$,
and another copy of $2$ before the end node $1$ which leads to
$Q=\seq{1,2,\ldots,2,1}$.
\item[(b)] If
2 is not a fixed node, we have two choices for the placement of 2: either next
to the start which leads to $Q=\seq{1,2,\ldots,1}$,  or next to the end which
leads to 
$Q'=\seq{1,\ldots,2,1}$.
\end{itemize}
After having placed item 2, we decide on the location of $3$, and so
on. At each step of the construction, the feasible sequence describes two
partially constructed feasible tours for the 2-TSP. 
The node placed in $Q$ on the final step  is node 1 which separates the nodes
in the two tours and serves as the last node in the first tour and as
the first node in the second tour.
 
 According to this description, the tours depicted in Fig.~\ref{fig:2pTSP} can be represented as the sequence
$Q=\seq{1,2,3,5,8,9,10,{\bf 1},8,7,6,5,4,3,1}$, with the set of fixed nodes $S=\{1,3,5,8\}$ and $s^*=8$. The reader who is familiar with the notion of \emph{ pyramidal} tours (see
e.g. Chapter 7 in \cite{GLS}) can recognise a pyramidal-like structure in the
sequence $Q$. A pyramidal tour consists of two parts: in the first part of the
tour the nodes are sorted in increasing order, in the
second part of the tour the nodes are sorted in decreasing order. 


To find an optimal sequence among all feasible sequences with the structure
described above, we use a
DP approach which is similar to the approach used for finding an optimal pyramidal tour. Assume that the nodes
$\{1,2,\ldots i,\ldots, j\}$, $i<j$,  have already got placed into the sequence
$Q$. Furthermore, assume that
node $i$ is currently the \emph{last} node in the first partially constructed
tour, and $j$ is currently the \emph{first} node in the second partially
constructed tour. Note that in this case the nodes $i+1,\ldots,j-1$ have
already been placed in the second tour.  

Let $V(i,j,m)$ denote the minimal length of a feasible subsequence 
that starts at node $i$, goes through the nodes $\{j+1,\ldots,n\}\cup\{1\}$ and
stops at node $j$ where $m$ denotes the number of nodes in
the second partially constructed tour in the sequence, i.e. the sub-path
$\seq{j,j-1,\ldots,i+1,\ldots,1}$.  
For the placement of $j+1$ we have to distinguish two cases.
If $j+1$ is a fixed node, it has to be placed in both tours, after $i$ and before $j$. If $j+1$ is not a fixed node, one has to decide whether to place
it  after $i$ or before $j$ depending on which option leads to a
shorter total length.

If there exists an item $s$ in the set of fixed items $S$ such that $i<s<j+1$,
it is not possible to find a \emph{feasible} subsequence from $i$ to $j$.
In this case we define $V(i,j,m):=\infty$.

Values $V(j,i,m)$ are defined similarly for the subsequences that start at $j$ and end at $i$. For the case when $j$ is a fixed node, we also need to define $V(j,j,m)$ as the length of the shortest feasible sequence from $j$ to $j$. 

It follows from the definitions above that the length of an optimal sequence $Q$, i.e.\ the total length of an optimal pair of tours, can be calculated as
\begin{equation}
\label{eq:vi3}
\begin{aligned}
c(Q)=
  \begin{cases}
  \min \{c_{12}+V(2,1,1),c_{21}+V(1,2,2)\}, 
\mbox{if}\ 2 \notin S;\\
 c_{12}+c_{21}+V(2,2,2),  \ \ \mbox{if}\ 2 \in S.
  \end{cases} 
  \end{aligned}
\end{equation}

The  $V$ values which involve node $n$ can be defined as shown below.

\begin{equation}
\label{eq:vi1}
\begin{aligned}
V(i,n,m)=
    c_{i1}+c_{1n},  \ \mbox{if}\ \lceil(n+|S|-p)/2\rceil\le m\le \lfloor(n+|S|+p)/2\rfloor, \\ 
V(n,i,m)=
    c_{n1}+c_{1i}, \ \mbox{if}\ \lceil(n+|S|-p)/2\rceil\le m\le \lfloor(n+|S|+p)/2\rfloor, \\
\qquad i=s^*,s^*+1,\ldots,n.\hspace{2cm}\\
\end{aligned}
\end{equation}

The above definitions ensure that balanced sequences that include all fixed nodes are considered.
It is easy to check that the remaining values $V(i,j,m)$ and $V(j,i,m)$ satisfy
the following DP recursions where we assume that undefined
values $V$  are set to infinity.

\begin{equation}
\label{eq:vi4}
\begin{aligned}
\hspace{-2ex}
V(i,j,m)=
\begin{cases}
  \min
  \begin{cases}
    c_{i,j+1}+V(j+1,j,m)\\
        \vspace{-0.4cm}
    \hspace{5.7cm}\mbox{if}\ j+1\notin S;\\
    c_{j+1,j} +V(i,j+1,m+1)\ \ \ \ 
     \end{cases}\\  
     c_{i,j+1}+c_{j+1,j}+V(j+1,j+1,m+1)\ \  \mbox{if}\ j+1 \in S,
\end{cases}\\
V(j,i,m)=
\begin{cases}
  \min
  \begin{cases}
    c_{j,j+1}+V(j+1,i,m)\\
    \vspace{-0.4cm}
    \hspace{5.7cm}\mbox{if}\ j+1\notin S;\\
    c_{j+1,i} +V(j,j+1,m+1)\ \ \ \  
     \end{cases}\\  
     c_{j,j+1}+c_{j+1,i}+V(j+1,j+1,m+1)\ \  \mbox{if}\ j+1 \in S,
\end{cases}\\
\qquad \forall\ i,j :\ \nexists\ s\in S \ \mbox{with}\   i<s<j+1.
\end{aligned}
\end{equation}
%

Summarizing this proves the main result of the section.
\begin{proposition}
The balanced 2-TSP on a Kalmanson matrix can be solved in $O(n^3)$  
time.
\end{proposition}

Note that the recursions (\ref{eq:vi3})--(\ref{eq:vi4}) could be further simplified
by taking into account the symmetricity of the distance matrix $C$. We believe
though that formulating the recursions as done above provides better insight
into the underlying structure.

\section{ Heuristic motivated by the Kalmanson solvable case} \label{sec:KS}
Given an instance of the 2-TSP with 
a symmetric $n\times n$ distance matrix $C$ and a subset $S$ of customers to be
visited twice, 
one can find a \emph{feasible} solution to the problem by using 
recursions (\ref{eq:vi3})--(\ref{eq:vi4}). The so obtained solution is
the best among
exponentially many feasible solutions and is optimal if $C$ is a Kalmanson
matrix. Obviously, for an arbitrary matrix $C$ the solution found is not necessarily optimal.
  Clearly, permuting rows and columns in matrix $C$, or in other words renumbering the cities in the problem, yields a  new exponential
neighbourhood searched in (\ref{eq:vi3})--(\ref{eq:vi4}), and hence a potentially new feasible solution. We describe here a heuristic that uses this idea for finding a solution for the 2-TSP.

One of the intuitive approaches to number the cities is to choose for this a heuristic solution to the TSP. To find a good TSP tour, we use a combination of the simplest and the best know TSP heuristics, nearest neighbour heuristic and the 2-opt heuristic (see e.g.\ Chapters 5 and 7 in \cite{GLS}). 
The tour found is used for renumbering the rows and columns in the distance matrix, and  recursions (\ref{eq:vi3})--(\ref{eq:vi4}) are then used to find a solution to the 2-TSP.

This heuristic is formally presented in the pseudo-code below.
We call it Kalmanson Sequence (KS-)Heuristic. The heuristic takes as input an $n\times n$ matrix $C$, subset $S$,  node {\bf start}, and returns a
sequence $Q(S)$ which is a heuristic solution to the 2-TSP with 
the distance matrix $C$. The points from set $S$ are visited in each of the two tours in $Q$. 

\begin{tabbing}
1111\=1111\=1111\=111\=111\kill
\bf KS-Heuristic for the balanced 2-TSP($n$, $C$, $S$, \bf{start}, $Q$)\\
{\bf \{}
Starting with node {\bf start} find the nearest neighbour tour on the set $\{1,\ldots,n\}$;\\
\hspace{4pt} Apply 2-opt to improve the tour; use obtained tour $\tau$ for renumbering the cities;\\
\hspace{6pt}Construct sequence $Q(S)$ with two 2-TSP tours \\
\>\> 
 by applying recursions (\ref{eq:vi3})--(\ref{eq:vi4}) to the permuted matrix 
$C_{\tau}$;\\
\hspace{4pt} Apply 2-opt to each of the tours in $Q(S)$; save the record;\\
\hspace{6pt}{\bf return} the record, i.e.\ the best $Q(S)$ found;\\
{\bf \}}
\end{tabbing}


We run computational experiments to test the performance of the 
KS-Heuristic. In the experiments we recorded the best solutions found by KS in $n$ runs (on the number of possible start  points in the nearest neighbour heuristic). 

The only published benchmark instances for the balanced 2-TSP are the  instances from Bassetto \& Mason \cite{BassettoM11}.
They considered a Euclidean version of the balanced 2-TSP.
A short summary of their approach  is as follows. 
First a TSP tour on the set of all customers, called a general tour (GT), is
constructed. The GT is used to obtain a partition of the set of customers
into two sub-tours. The initial partition of GT into two sub-tours is
improved by applying decision rules motivated by geometry (e.g. solutions get improved by removing crossing edges).
For each of the sub-tours, an optimal TSP tour is constructed by applying an exact TSP algorithm (see \cite{ABCC}, Chapter 16,  and \cite{Concorde}). 
The authors also use visualisation and human intervention for improvements of
the solutions found by the algorithms.

The set of benchmark instances from \cite{BassettoM11} contains 60 randomly
generated instances with 48 customers each. The set of these instances is divided
into three subsets with a different number of customers to be visited in both periods: $8$, $12$, and $24$ customers.   
We use these instances to test our KS-Heuristic.
The summary of the computational experiments is presented in Table \ref{table:summaryHeuristics}.
We compare our results with both types of solutions presented in
\cite{BassettoM11}: with the best solutions found by a computer, which we refer to as
``PC" solutions, and  with the solutions obtained after visualisation and additional human intervention
referred to  as ``PC/h". 
  For each set of 20 instances, Table~\ref{table:summaryHeuristics} shows the
  mean of the percentage above the best solution, the percentages for the
  ``best" and ``worst" instance, and the number of instances where we improved
the previously best solution, found solutions with the same value, or failed to
improve the best known solution (cf. the rows \#($<,=,>$) in the tables). A 
negative percentage means that the best known solution was improved. Results for individual test instances can be found in  Appendix A.

For our experiments we used a laptop with an 
Intel(R) i5-8350U 1.70 GHz CPU and 16 GB of RAM. The alogorithms were
implemented in
C++ and compiled with MinGW-w64. Computational
times for the instances were less than 30 milliseconds (\cite{BassettoM11}
mentioned times of a ``few seconds'' in their experiments).

As Table~\ref{table:summaryHeuristics} shows,  the suggested 
heuristic improves solutions for  22 out of 60 instances tested in
\cite{BassettoM11}, with the quality of other solutions not far away from 1\% accuracy. 
For the case with 8 fixed points our results are on overall better than the results in  \cite{BassettoM11}.
This backs up our belief that polynomially solvable cases
of hard optimization problems can  provide the basis for the construction of
heuristics with a performance  which is competitive with other approaches.

For the solutions  obtained  by visualisation of the computer solutions and
follow-up human intervention our heuristic improved  6 out of 60
results. Still the average quality of our solutions which are obtained in milliseconds is within 2\% from the previously found solutions. In
the next section we describe a new heuristic which is able to further
improve many of the published solutions.

The DP recursions (\ref{eq:vi3})--(\ref{eq:vi4}) make use of the structure ``one sequence - two tours" and simultaneously optimise an allocation of nodes to tours and the sequence of the nodes in the tours.  We generalise this approach and develop a heuristic which is based on the same methodology and uses the well-known Held and Karp \cite{HeldKarp} DP recursions for the TSP with arbitrary distance matrices.

\begin{table}
\caption{Summary of computational experiments for the KS-Heuristic.}\label{table:summaryHeuristics}
\begin{center}
\small
\begin{tabular}{||c||l||c|c||}
\hline 
\cline{2-4} 
& 
& PC & PC/h  
\\ 
\hline 
&Mean \%   &-0.56  & 1.18  \\ 
8 nodes& Best \% & -6.11  & -1.18  \\ 
visited twice&Worst \%  &3.7  & 4.92   \\ 
& \# (\scalebox{0.8}{$<,=,>$}) & (15,0,5) & (3,0,17) \\ 
\hline 
&Mean \%   & 0.72  & 1.78  \\ 
16 nodes& Best \% & -1.74  & -1.54  \\ 
visited twice&Worst \%  &4.21  & 6.31   \\ 
& \# (\scalebox{0.8}{$<,=,>$}) & (6,0,14) & (3,0,17) \\ 
\hline 
&Mean \%  & 1.25  & 1.88  \\ 
24 nodes& Best \%   & -0.04  & 0.32   \\ 
visited twice&Worst \%   & 3.16  & 4.65   \\ 
&\# (\scalebox{0.8}{$<,=,>$}) & (1,0,19) & (0,0,20) \\ 
\hline 
Total&  & (22,0,38) & (6,0,54) \\ 
\hline 
\end{tabular} 
\end{center}
\end{table}%

\section{``One sequence - two tours" heuristic for the 2-TSP}
\label{sec:framework}
\subsection{  Formulation of the 2-TSP as a vehicle routing problem}
In this section we model the 2-TSP within the framework of a vehicle
routing problem (VRP) with two vehicles (2-VRP).
We believe that the 2-VRP model is better suited for presenting our ideas.

Consider the 2-TSP with $n$ nodes with $|S|$ of them to be visited twice. Let one of the nodes from set $S$ be a depot and renumber the rest of the nodes as 1,2,\ldots,$n-1$. Let $F_1\subset\{1,\ldots,n-1\}$ be the set of nodes to be visited twice. Define $F_2:=F_1$ and then renumber the nodes in $F_2$ as $\{n,n+1,\ldots,N:=n+|S|-2\}$. In the context of the 2-VRP, set $F_1$ is the set of fixed customers to be visited by vehicle one, $F_2$ is the set of fixed customers to be visited by vehicle two.

Let $W_1$ and $W_2$ be the capacities of vehicles one and two. If  each customer $i$ has a unit demand $w(i)=1$, then the balancing constraint requiring the same number of customers visited by each vehicle can be met by  setting the capacities
as $W_1=W_2=\lfloor(n-1+|S|+p)/2\rfloor$, with $p=1$. 

In this manner we end up with a 2-VRP instance with one depot and the set of customers $\{1,2,\ldots,N\}$.  Set $F_1$ is the set of customers to be always allocated to the vehicle one, and  $F_2$ is the set of customers to be always allocated to the vehicle two ($F_1\cap F_2=\emptyset$). Notice, in real life applications one often needs to allocate customers to particular vehicles, e.g.\ due to a restricted site access or drivers' preferences.

\subsection{ Vehicle routing problem with customers as subsets}

Our approach to solving the 2-VRP uses the well-known Held \& Karp \cite{HeldKarp}  DP algorithm for the TSP. This is an exact algorithm  and naturally can solve only small size instances. To deal with large size instances, we will repeatedly combine several customers into a single customer. With a new customer we will associate a new demand which is the total demand of the customers combined. We also need to consider an extra travel distance while visiting a new customer. It may be helpful to use the following interpretation. 
 
We assume that each customer $i$ is located in an estate with only two entry points.
To distinguish between two entry points to the estate, we refer to  one of the
entry points as the {\em left end} and denote it by $L(i)$. The other of the
entries  is  referred to as the {\em right end}, and is denoted by $R(i)$. 
Representing the locations in the described way makes the task of data
aggregation more intuitive. Assume that there is a collection of customers
which are located in the same estate. If it is possible for these customers to be serviced by
one vehicle, then all these customers can be substituted by one new customer whose demand
equals the total demand of the customers. The distance from $L(i)$ to $R(i)$
(and from $R(i)$ to $L(i)$, since the distances are symmetric), denoted as
$l(i)$, can be calculated as the length of the shortest path visiting all customers in the estate. 

Summarising, each customer in our model will have four attributes:
the left end $L(i)$, the right end $R(i)$, the distance $l(i)$ of travelling between the left and the right ends, and the demand $w(i)$.

It will be convenient to assume that one vehicle travels from depot $d^1$ to depot $0$, and another vehicle travels from depot $0$ to depot $d^2$. Depot $0$ will be considered as a special customer with the set of attributes
$\{L(0)$, $R(0)$, $l(0):=0$, $ w(0):=0\}$. 
The reason for introducing the customer $0$ is to use it for separating customers and allocating them to the two different vehicles. This restricts positions for placing customer $0$: the total demand of customers visited by each vehicle should not exceed the vehicle capacities. This constraint can be easily dealt with in the DP recursions.

The main idea of our approach is to view two routes for the two vehicles as a \emph{single two-vehicle} route (``one sequence - two routes"), and  then to apply an algorithm for the TSP to find an optimal route. In the literature on the VRP the term \emph{route} is used widely, therefore from now on we use  both terms, \emph {route} and \emph{tour}.

\subsection{\bf Dynamic programming recursions}
In this section we adapt the Held \& Karp \cite{HeldKarp} DP recursions for finding the shortest path from $d^1$ to $d^2 $ through the set of customers  $U=\{0\}\cup\{1,2,\ldots,N\}$. In this path the customers from set $F_1$ are placed before customer 0, the customers from set $F_2$ are placed after customer 0. Moreover, the total demand of customers placed before customer 0 does not exceed capacity $W_1$, the total demand of customers placed after customer 0 does not exceed capacity $W_2$.

Let $J$ be a subset of customers not containing $i$, so $J\subseteq U\setminus \{i\}$.
  Let $\VL[i,J]$ be the minimum length of an optimal 2-vehicle route among all routes that start visiting customer $i$ from the left end, then visiting all the customers in set $J$, and stopping in depot $d^2$. Similarly, define $\VR[i,J]$ to be the length of the optimal route that starts visiting customer $i$ from the right end. 
  Although we are dealing with symmetric distance matrices here,  it is important while sequencing the customers to distinguish whether we pass a customer from the left or from the right end. The optimal length of the 2-vehicle route can be calculated as 
\begin{equation}
\label{eq:dp0}
V=min_{i\in U\setminus \{\{0\}\cup F_2\}}\{c_{d^1,L(i)}+\VL[i,U\setminus \{i\}],c_{d^1,R(i)}+\VR[i,U\setminus \{i\}]\}.
\end{equation}

In the formula above, $i$ is the first customer in the path. We assume that the total demand of customers is greater than the capacity of each of the vehicles, so item $0$ cannot be the first customer. Items from set $F_2$ are also excluded from the  consideration in this step.

Values $\VL[i,J]$ and $\VR[i,J]$ for all items $i$ and subsets $J$, $J\subset \{0\}\cup\{1,2,\ldots,N\}$, are calculated as shown in the recursions below. In the preprocessing step we calculate the total demand $w(J)$ for subsets $J$ and define initially all values $\VL[i,J]$ and $\VR[i,J]$ to be  infinity.

\begin{equation}
\label{eq:dp1}
\begin{aligned}
\VL[i,J]_{i\neq 0}& =
\begin{cases}
\left.
\min_{j\in J\setminus F_2}
\begin{cases}
				l(i)+c_{R(i),L(j)}+\VL[j,J\setminus\{j\}] \\
				l(i)+c_{R(i),R(j)}+\VR[j,J\setminus\{j\}]
\end{cases}
\hspace{-2.0ex}\right\}
&\hspace{-2.0ex}\mbox{if }0\in J \\
\\
\left.
\min_{j\in  J\setminus F_1}
\begin{cases}
				l(i)+c_{R(i),L(j)}+\VL[j,J\setminus\{j\}]\\
				l(i)+c_{R(i),R(j)}+\VR[j,J\setminus\{j\}]
\end{cases}
\hspace{-2.0ex}\right\}
&\hspace{-2.0ex}\mbox{if }0\notin J, w(\{i\}\cup J)\le W_2
\end{cases}
\end{aligned}
\end{equation}

\begin{equation}
\begin{aligned}
\VR[i,J]_{i\neq 0}& =
\begin{cases}
\left.
\min_{j\in J\setminus F_2}
\begin{cases}
				l(i)+c_{L(i),L(j)}+\VL[j,J\setminus\{j\}] \\
				l(i)+c_{L(i),R(j)}+\VR[j,J\setminus\{j\}]
\end{cases}
\hspace{-2.0ex}\right\}
&\hspace{-2.0ex}\mbox{if }0\in J, \\\\
\left.
\min_{j\in J\setminus F_1}
\begin{cases}
				l(i)+c_{L(i),L(j)}+\VL[j,J\setminus\{j\}]\\
				l(i)+c_{L(i),R(j)}+\VR[j,J\setminus\{j\}]
\end{cases}
\hspace{-2.0ex}\right\}
&\hspace{-2.0ex}\mbox{if }0\notin J, w(\{i\}\cup J)\le W_2
\end{cases}
\end{aligned}
\end{equation}

\begin{equation}
\begin{aligned}
\VL[0,J]& =
\left.
\min_{j\in J\setminus F_1}
\begin{cases}
				c_{R(0),L(j)}+\VL[j,J\setminus\{j\}\\
				c_{R(0),R(j)}+\VR[j,J\setminus\{j\}
\end{cases}
\hspace{-2.0ex}\right\}
&\mbox{if }
w(U\setminus  J)\le W_1,
w(J)\le W_2.
\end{aligned}
\end{equation}

The boundary conditions are:
\begin{equation}
\label{eq:dpB}
\begin{aligned}
\VL[i,\emptyset]_{i\notin \{0\}\cup F_1}& =
				l(i)+c_{R(i),d^2},
\  \  
\VR[i,\emptyset]_{i\notin \{0\}\cup F_1}& =
				l(i)+c_{L(i),d^2}.
\end{aligned}
\end{equation}

 In the recursions above, if subset $J$ does not contain item 0, we eliminate from the calculations the items which are in $F_1$. If subset $J$  contains item 0, we eliminate the items which are in $F_2$. It ensures that the fixed items are properly positioned.  Since there are only two vehicles, the capacity constraint is also easily checked without any additional calculations.

For instances with small sizes ($N\le 20$ in our experiments) the system
(\ref{eq:dp0})--(\ref{eq:dpB}) was solved within seconds. However when the
number of customers increases, the computations turn out to be too time
consuming which is hardly surprising. To make the DP approach
practical for larger instances as well, we suggest the heuristic approach
described below.

\subsection{\bf Aggregation strategy and local search}
\label{sec:aggregation}

\begin{figure}[tb]
\unitlength=0.8cm
\begin{center}
\begin{picture}(11.5,6.5)
\put(-2.5,0)
{
\begin{picture}(10.5,6.)
\includegraphics[scale=0.9]{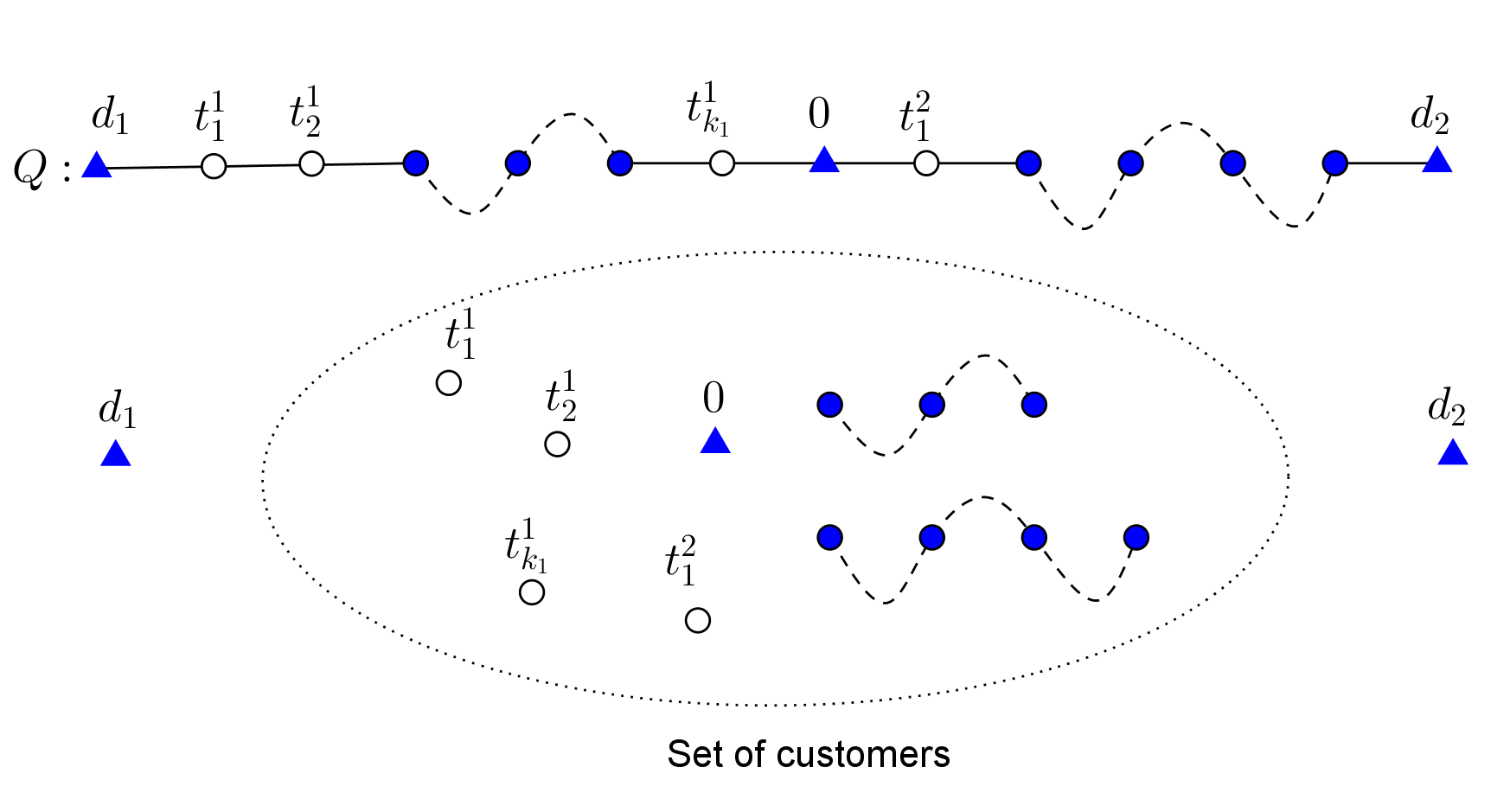}
\end{picture} 
}
\end{picture}
\end{center}
\caption{ Illustration of the ``sliding subset" heuristic: first step of disassembling; each subset contains $2$ customers, $S_1=\{t_1^1,t_2^1\}$, $S_2=\{t_{k_1}^1,t_1^2\}$.}
\label{fig:disassembling1}
\end{figure}

Consider the 2-VRP with $N$ customers, where $N$ is big enough to make the DP recursions computationally intractable. 
We start with a simple heuristic to find a feasible solution.  We then ``disassemble"  this solution into a small number of sub-paths, and  
represent each sub-path as a customer. 
Recalculation of the attributes associated with the new customer $i$ is
straightforward. The left end $L(i)$ for the new customer is the first node
in the sub-path; the right end $R(i)$ of the new customer is the last node in the
sub-path; the demand $w(i)$  is the total demand of customers in the
sub-path; value $l(i)$ is the length of the sub-path.

We then apply recursions (\ref{eq:dp0})--(\ref{eq:dpB}) and find an optimal
solution to the small-size 2-VRP with the new set of customers.   
The optimal solution found can now be transformed into a solution for the initial 2-VRP by replacing the aggregated customers with the original sub-paths.
Obviously, the so obtained solution to the initial 2-VRP instance
will in general not be optimal. 
Note however that the obtained solution can be  ``disassembled" again to
get a new small-size instance. 
This process is repeated until  a certain stopping criterion is satisfied.

\smallskip
We suggest the following approach for disassembling, which we call the sliding subset method. 
Assume that we have an initial solution  to the 2-VRP which we present as a two vehicle route: $Q=\seq{d_1,t^1_1,t^1_2,\ldots,t^1_{k_1},0,t^2_1,t^2_2,\ldots,t^2_{k_2},d_2}$. Here the route of vehicle 1 is 
$\seq{d_1,t^1_1,t^1_2,\ldots,t^1_{k_1},0}$ and the route of vehicle 2 is $\seq{0,t^2_1,t^2_2,\ldots,t^2_{k_2},d_2}$. Fig.~\ref{fig:disassembling1} illustrates the concept of sliding subset for an instance with $k_1=6$ and $k_2=5$. 

We disassemble the initial solution into a new set of customers as follows.
First, we choose customer $0$ as a customer in the new 2-VRP instance, and delete it 
from $Q$.
Let $s$ be a small integer constant, a parameter of the algorithm ($s=2$ in Fig.~\ref{fig:disassembling1}). Choose two subsets of customers $S_1$ and $S_2$  containing $s$ items each. The customers in each subset are at consecutive positions in $Q$. Subset $S_1$ will always contain at least one customer from the route of vehicle 1, and $S_2$ contains at least one customer from the route of vehicle 2. By defining subsets in this way we try to always have a possibility to move the customers between the vehicles (unless the customers are fixed).

In the first disassembling phase  define  $S_1=\{t^1_1,\ldots,t^1_s\}$, and
$S_2=\{t^1_{k_1-s+2},\ldots,t^1_{k_1},t^2_1\}$. 
Delete
the nodes in $S_1$ and $S_2$ from $Q$ and add them to the set of 
customers in the new 2-VRP instance. Depot $d_2$ stays as depot in the new problem; 
delete it from $Q$. 
The sub-paths which are left in $Q$ after the deletions will be replaced by  
aggregated customers and added to the new 2-VRP instance. 
These customers are: depot $d_1$ (as a sub-path with a single node), the sub-path $\seq{t^1_{s+1},\ldots,t^1_{k_1-s+1}}$ and the sub-path
$\seq{t^2_2,\ldots,t^2_{k_2}}$.

We \emph{redefine} the subset $S_2$ by deleting $l$ first elements from $S_2$ and adding $l$ new elements; $l$ is a parameter, to which we refer as {\em step}.  We repeat the process until we reach the end of the route. 
When subsets $S_2$ of the described type have been enumerated, we change set $S_1$ (with the step $l$) and redefine set $S_2$ to follow set $S_1$ in a manner as described above.

If the solution to the initial 2-VRP instance is improved, it is saved as the
record and the process of disassembling is applied to the new solution. The
procedure stops when all feasible subsets $S_1$ and $S_2$ got enumerated and no
improvements were found.

To keep the size of all small 2-VRP instances the same, we slightly modified
the disassembling step used in our implementation (for further details see
Appendix B).

\section{Computational experiments}
\label{sec:emp}

For the sliding subset heuristic with 
the size of sliding subsets $s$ and  step $l$, we use the notation
$H_{s,l}$. In computational experiments reported in this section we test the
influence of parameters $s$ and $l$ on the quality of solutions
obtained as well as the importance of initial solutions.
The algorithm used in the experiments is formally described below.

\begin{tabbing} 
1111\=1111\=1111\=111\=111\kill
\bf One Sequence - Two Tours  Heuristic for the 2-VRP($n$, $C$, $s$, $l$)\\
{\bf \{}
\>\emph{ \bf until}  (stop criteria is satisfied) { \bf do} \\
\>\emph{ \bf \{ }
Generate next feasible sequence $Q$ containing two tours of the  2-TSP; \\
\>\>\emph{\bf while} improvement found { \bf do} \\
\>\>\emph{ \bf \{  }
 \hspace{1pt} 
 Apply recursions (\ref{eq:dp1})--(\ref{eq:dpB}) implemented in $H_{s,l}$ to improve the tours in $Q$;\\
\>\>\>Apply 2-opt to each of the two tours; save the record;\\
\>\>\>Swap tours in $Q$ to apply (\ref{eq:dp1})--(\ref{eq:dpB}) to the new ordering in the next iteration;\\
\>\>\emph{ \bf \} }\\
\>\emph{ \bf \} }\\
\>{\bf return} record, i.e.\ the best $Q$ found;\\
{\bf \}}
\end{tabbing}

The stop criteria in the algorithm above will be defined as the number of iterations or as the total time limit for all iterations.

On each step in the improvement loop of the algorithm, recursions (\ref{eq:dp1})--(\ref{eq:dpB}) are used repeatedly in the disassemble-assemble steps in heuristic $H_{s,l}$ as described in Section~\ref{sec:aggregation}. The outcomes of these steps can be different, if the order of tours in sequence $Q$ is changed. Hence the swap step in the algorithm.


Two different approaches are used to generate feasible sequences $Q$. In the first set of experiments, the KS-heuristic (see Section~\ref{sec:KS}) is used to generate sequences $Q$. The number of \emph{start} points in the KS-Heuristic is  the stop criteria in this case.
The number of start points is bounded by the number of nodes in the problem. Moreover, for some of the start points the combination of the nearest neighbour and 2-opt heuristics may produce the same solutions. 
On the other hand, there is an exponentially many  feasible solutions partitioning the initial set of points into two routes in the sequence $Q$. Motivated by these observations we also considered a random partitioning of nodes in two routes as an alternative procedure for generating sequences $Q$. In this case the number of possible iterations is practically unrestricted. It can be restricted  by a chosen constant or by the bound on the total time the algorithm is allowed to run.

Several sliding subset heuristics have been chosen for the experiments.
To manage the computational time, different numbers of initial solutions have been chosen for different heuristics. The heuristics are:   $H_{4,2}$ (run with 48 initial solutions), $H_{5,3}$ (36 initial solutions), $H_{6,4}$ (24 initial solutions), and $H_{7,5}$ (12 initial solutions). 
The set of 60 benchmark  problems from \cite{BassettoM11} have been used.

Table~\ref{table:summary24_36} summarises results of the first set of experiments where the initial solutions have been obtained by the KS-heuristic. Notation $[KS+H_{4,2}]\times 48$ means that the corresponding columns contain the summary of results where the KS heuristic was followed by $H_{4,2}$; 48 is the number of generated initial solutions. The other notations in the table are the same as the notations used in Table~\ref{table:summaryHeuristics} in Section~\ref{sec:KS}.

In the second set of computational experiments we tested the importance of good initial solutions. Instead of KS-Heuristic for generating initial solutions the following Random Partition (RP) heuristic is used.
\needspace{8\baselineskip}
\begin{tabbing}
1111\=1111\=1111\=111\=111\kill
\bf RP-Heuristic to get a random feasible 2-TSP solution($n$, $C$, $S$, $Q$)\\
{\bf \{}
\>Create a copy of set $S\setminus\{1\}$; add the copies to the set $\{1,\ldots,n\}$;\\
\>
Use C++ function to \emph{random\_shuffle}  the set of indices $\{2,\ldots,n+|S|-1\}$;\\
\>Transform the obtained permutation into sequence $Q$ with two tours:\\
\>\>points from $S$ are left in the first tour, and their copies are in the second; \\
\> Apply 2-opt to each of the tours in $Q$; \\
{\bf \}}
\end{tabbing}

Using the info from Table~\ref{table:summary24_36} we aimed to reach the same computational time as was spent on instances in the first set of experiments. Therefore the stop criteria in this case was chosen to be the computational time spent in the first set of experiments.
The mean time achieved in experiments could still be slightly different from the target, so we keep this info in the table. 
The results of the experiments are summarised in Table~\ref{table:summary12_24}.

\begin{table}[htbp]
  \centering
  \caption{Test results with initial solutions obtained by Kalmanson Sequence heuristic}
  \small
\begin{tabular}{|c|l|c|c|c|c|c|c|c|c|}
\hline
\cline{3-10}\multicolumn{1}{|c}{} &     & \multicolumn{2}{c|}{[KS+$H_{4,2}]\times 48$} & \multicolumn{2}{c|}{[KS+$H_{5,3}]\times 36$} & \multicolumn{2}{c|}{[KS+$H_{6,4}]\times 24$} & \multicolumn{2}{c|}{[KS+$H_{7,5})]\times 12$} \bigstrut\\
\cline{3-10}\multicolumn{1}{|c}{} &     & PC  & PC/h & PC  & PC/h & PC  & PC/h & PC  & PC/h \bigstrut\\
\hline
    &     & \multicolumn{2}{c|}{$t_m\approx 29s$} & \multicolumn{2}{c|}{$t_m\approx 70s$} & \multicolumn{2}{c|}{$t_m\approx 81s$} & \multicolumn{2}{c|}{$t_m\approx 166s$} \bigstrut\\
\cline{3-10} 8 nodes 
 & Mean \% & $-\bf{ 2.37}$ & $-{\bf 0.66}$ & $-$2.41 & $-$0.70 & $-$2.44 & $-0.73$ & $-${\bf 2.59} & $-${\bf 0.79} \bigstrut[t]\\
visited  & Best \scalebox{0.85}{\%} & $-$7.31 & $-$1.89 &$-$7.31 & $-$2.65 & $-$7.31 & $-$1.84 & $-$7.31 & $-$2.03 \\
twice & Worst \scalebox{0.85}{\%} & $-$0.36 & +0.51 & $-$0.36 &+0.87 & $-$0.24 & +0.94 & $-$0.39 & +0.11 \\
    & \scalebox{0.8}{\#($<,=,>$)} & \scalebox{0.85}{(20,0,0)} & \scalebox{0.85}{(14,2,4)} & \scalebox{0.85}{(20,0,0)} & \scalebox{0.85}{(15,3,2)} & \scalebox{0.85}{(20,0,0)} & \scalebox{0.85}{(16,2,2)} & \scalebox{0.85}{(20,0,0)} & \scalebox{0.85}{(18,0,2)} \\
    \hline
    &     & \multicolumn{2}{c|}{$t_m\approx 30s$} & \multicolumn{2}{c|}{$t_m\approx 83s$} & \multicolumn{2}{c|}{$t_m\approx 110s$} & \multicolumn{2}{c|}{$t_m\approx 206s$} \bigstrut[b]\\
\cline{3-10}16 nodes & Mean \scalebox{0.85}{\%} & $-$1.56  & $-$0.52  & $-$1.80  & $-$0.77  & $-$1.83  & $-$0.80  & $-${\bf 1.91}  & $-${\bf 0.88}      \bigstrut[t]\\
visited  & Best \scalebox{0.85}{\%}& $-$3.77  & $-$1.84  & $-$3.93 & $-$2.59   & $-$3.89  & $-$2.41  & $-$3.89  & $-$2.24  \\
twice & Worst \scalebox{0.85}{\%}& +0.22  & +1.00  & $-$0.10  & +1.23 & $-$0.30 & +0.94 & $-$0.49  & +0.70\\
    & \scalebox{0.8}{\#($<,=,>$)} & \scalebox{0.85}{(19,0,1)} & \scalebox{0.85}{(15,0,5)} & \scalebox{0.85}{(20,0,0)} & \scalebox{0.85}{(18,0,2)} & \scalebox{0.85}{(20,0,0)} & \scalebox{0.85}{(18,0,2)} & \scalebox{0.85}{(20,0,0)} & \scalebox{0.85}{(17,0,3)} \\
    \hline
    &     & \multicolumn{2}{c|}{$t_m\approx 32s$} & \multicolumn{2}{c|}{$t_m\approx 91s$} & \multicolumn{2}{c|}{$t_m\approx 90s$} & \multicolumn{2}{c|}{$t_m\approx 230s$} \bigstrut[b]\\
\cline{3-10}24 nodes & Mean \scalebox{0.85}{\%}& $-$0.44  & $+$0.18  & $-$0.97  & $-$0.35 & $-$0.98  & $-$0.36 & $-$0.97  & $-$0.35\bigstrut[t]\\
visited  & Best \scalebox{0.85}{\%} & $-$2.39  & $-$0.64  & $-$3.80  & $-$1.41  & $-$4.44  & $-$1.49  & $-$3.70  & $-$1.43  \\
twice & Worst \scalebox{0.85}{\%}  & +0.49  & +0.88  & +0.32  & +0.80  & +0.62  & +1.30  & +0.37  & +0.72 \\
    & \scalebox{0.8}{\#($<,=,>$)} & \scalebox{0.85}{(16.0,4)} & \scalebox{0.85}{(10,0,10)} & \scalebox{0.85}{(19,0,1)} & \scalebox{0.85}{(14,1,5)} & \scalebox{0.85}{(18,0,2)} & \scalebox{0.85}{(14,1,5)} & \scalebox{0.85}{(18,0,2)} & \scalebox{0.85}{(12,2,6)} \\
    \hline
    \hline
   Total     & \scalebox{0.8}{\#($<,=,>$)} & \scalebox{0.85}{(55,0,5)} & \scalebox{0.85}{(39,2,19)} & \scalebox{0.85}{(59,0,1)} & \scalebox{0.85}{(47,4,17)} & \scalebox{0.85}{(58,0,2)} & \scalebox{0.85}{(48,3,9)} & \scalebox{0.85}{(58,0,2)} & \scalebox{0.85}{(47,2,11)} \\
        \hline  
\end{tabular}%
  \label{table:summary24_36}%
  \hspace{1cm} 
  \centering
  \caption{Test results with initial solutions obtained by Random Partition heuristic}
  \small
\begin{tabular}{|c|l|c|c|c|c|c|c|c|c|}
\hline
\cline{3-10}\multicolumn{1}{|c}{} &     & \multicolumn{2}{c|}{[RP+$H_{4,2}] $} & \multicolumn{2}{c|}{[RP+$H_{5,3}]$} & \multicolumn{2}{c|}{[RP+$H_{6,4}]$} & \multicolumn{2}{c|}{[RP+$H_{7,5}]$} \bigstrut\\
\cline{3-10}\multicolumn{1}{|c}{} &     & PC  & PC/h & PC  & PC/h & PC  & PC/h & PC  & PC/h \bigstrut\\
\hline
    &     & \multicolumn{2}{c|}{$t_m\approx 29s$} & \multicolumn{2}{c|}{$t_m\approx 70s$} & \multicolumn{2}{c|}{$t_m\approx 81s$} & \multicolumn{2}{c|}{$t_m\approx 160s$} \bigstrut\\
\cline{3-10} 8 nodes & Mean \scalebox{0.85}{\%} & $-$2.20 & $-$0.49 & $-${\bf 2.45} & $-${\bf 0.74} & $-${\bf 2.56} & $-${\bf 0.82} & $-$2.52 & $-$0.81 \bigstrut[t]\\
visited  & Best \scalebox{0.85}{\%} & $-$6.90 & $-$2.51 & $-$7.19 & $-$2.65 & $-$6.85 & $-$2.44 & $-$7.31 & $-$2.77 \\
twice & Worst \scalebox{0.85}{\%} & +1.71 & +1.71 & +1.60 & +1.60 & +0.28 & +0.39 & +0.28 & +0.75 \\
    & \scalebox{0.8}{\#($<,=,>$)} & \scalebox{0.85}{(19,0,1)} & \scalebox{0.85}{(11,1,8)} & \scalebox{0.85}{(19,0,1)} & \scalebox{0.85}{(14,2,4)} & \scalebox{0.85}{(19,0,1)} & \scalebox{0.85}{(15,1,4)} & \scalebox{0.85}{(19,0,1)} & \scalebox{0.85}{(13,0,7)} \\
    \hline
    &     & \multicolumn{2}{c|}{$t_m\approx 29s$} & \multicolumn{2}{c|}{$t_m\approx 83s$} & \multicolumn{2}{c|}{$t_m\approx 109s$} & \multicolumn{2}{c|}{$t_m\approx 205s$} \bigstrut[b]\\
\cline{3-10}16 nodes & Mean \scalebox{0.85}{\%} & $-${\bf 1.75} & $-${\bf 0.71} & $-${\bf 1.88} & $-${\bf 0.85}& $-${\bf 2.09} & $-${\bf 1.07} & $-$1.90 & $-$0.87 \bigstrut[t]\\
visited  & Best \scalebox{0.85}{\%} & $-$3.75 & $-$2.41 & $-$3.75 & $-$2.41 & $-$3.75 & $-$2.41 & $-$3.86 & $-$2.52 \\
twice & Worst \scalebox{0.85}{\%} & $-$0.31 & +0.41 & $-$0.07 & +0.72 & $-$0.52 & +0.19 & $-$0.07 & +2.09 \\
    & \scalebox{0.8}{\#($<,=,>$)} & \scalebox{0.85}{(20,0,0)} & \scalebox{0.85}{(16,0,4)} & \scalebox{0.85}{(20,0,0)} & \scalebox{0.85}{(17,0,3)} & \scalebox{0.85}{(20,0,0)} & \scalebox{0.85}{(17,0,3)} & \scalebox{0.85}{(20,0,0)} & \scalebox{0.85}{(16,1,3)} \\
    \hline
    &     & \multicolumn{2}{c|}{$t_m\approx 31s$} & \multicolumn{2}{c|}{$t_m\approx 91s$} & \multicolumn{2}{c|}{$t_m\approx 90s$} & \multicolumn{2}{c|}{$t_m\approx 229s$} \bigstrut[b]\\
\cline{3-10}24 nodes & Mean \scalebox{0.85}{\%} & $-${\bf 0.93} & $-${\bf 0.31} & $-${\bf 1.08} & $-${\bf 0.46} & $-${\bf 1.02} & $-${\bf 0.40} & $-${\bf 1.03} & $-${\bf 0.41} \bigstrut[t]\\
visited  & Best \scalebox{0.85}{\%} & $-$3.34 & $-$1.50 & $-$3.67 & $-$1.11 & $-$3.48 & $-$1.52 & $-$3.98 & $-$1.52 \\
twice & Worst \scalebox{0.85}{\%} & +0.64 & +0.94 & +0.76 & +0.76 & +0.76 & +0.76 & +0.76 & +0.76 \\
    & \scalebox{0.8}{\#($<,=,>$)} & \scalebox{0.85}{(16,0,4)} & \scalebox{0.85}{(12,1,7)} & \scalebox{0.85}{(18,0,2)} & \scalebox{0.85}{(15,1,4)} & \scalebox{0.85}{(19,0,1)} & \scalebox{0.85}{(15,0,5)} & \scalebox{0.85}{(19,0,1)} & \scalebox{0.85}{(13,2,5)} \\
    \hline
    \hline
   Total     & \scalebox{0.8}{\#($<,=,>$)} & \scalebox{0.85}{(55,0,5)} & \scalebox{0.85}{(39,2,19)} & \scalebox{0.85}{(57,0,3)} & \scalebox{0.85}{(46,3,11)} & \scalebox{0.85}{(58,0,2)} & \scalebox{0.85}{(47,1,12)} & \scalebox{0.85}{(58,0,2)} & \scalebox{0.85}{(46,2,12)} \\
        \hline
    
\end{tabular}%
  \label{table:summary12_24}%
\end{table}%

\begin{figure}
\unitlength=1cm
\begin{center}
{\begin{picture}(15.5,8.)
\hspace{-0.7cm}
\includegraphics[scale=0.22]{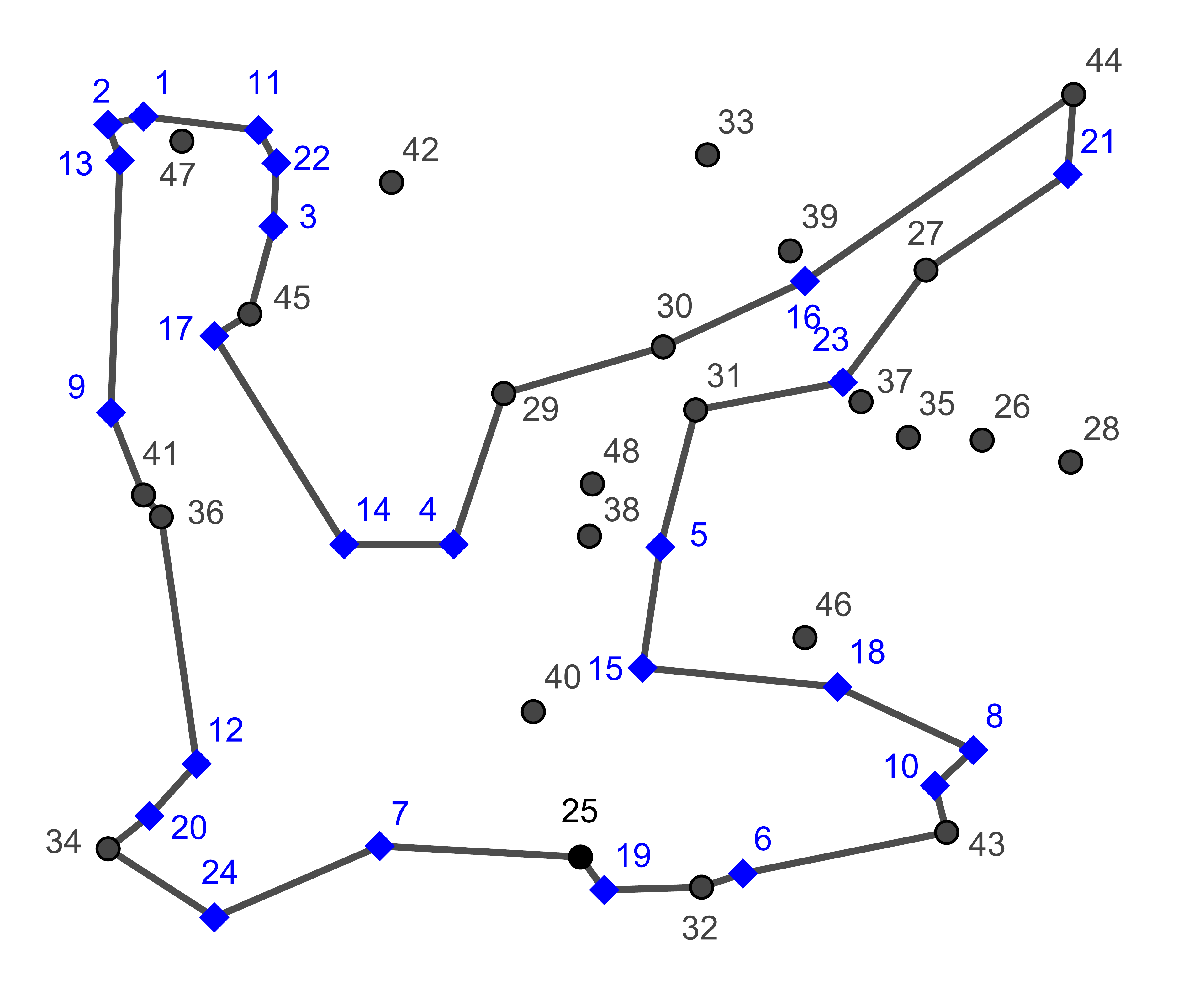}
\includegraphics[scale=0.22]{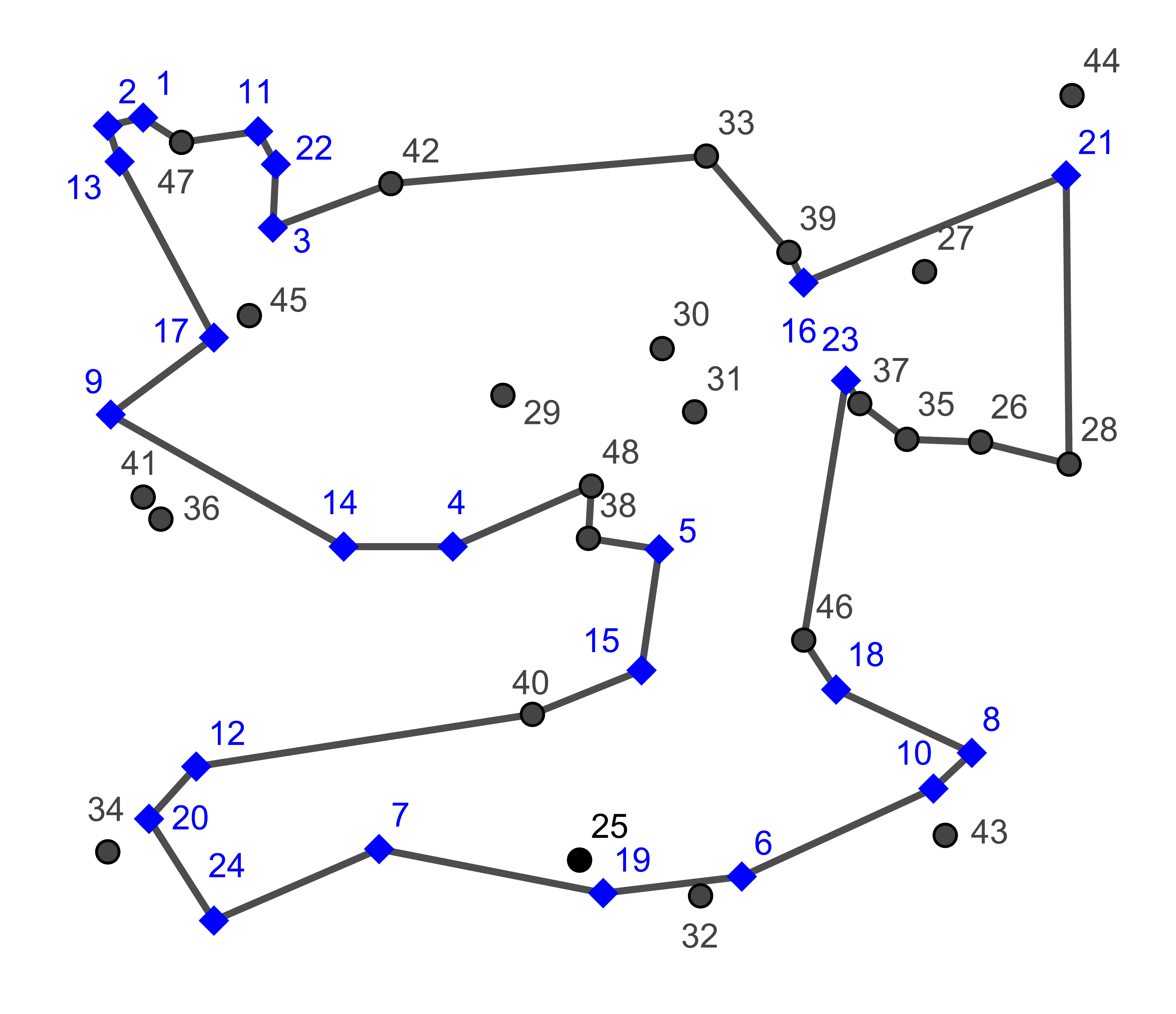}
\end{picture} }
\end{center}
\caption{Instance 52: Solution with the length 32201 (0.24\% improvement) 
 }
\label{fig:Instance52}
\end{figure}

Analysing the results of computational experiments and the information presented in Tables~\ref{table:summary24_36} and \ref{table:summary12_24} (see also the table in Appendix A) we can draw the following conclusions.
\begin{itemize}
\item By changing parameters $s$ and $l$ in the heuristics $H_{s,l}$ it is possible to trade off between the computational time and accuracy of algorithms. Increasing the number of iterations is another possibility to influence the quality of the solutions obtained;
\item Surprisingly, algorithms with random initial partitions performed better than the algorithms with the initial solutions obtained with the KS-heuristic (better results in the tables are highlighted with the bold font). It means that the 2-VRP heuristic is robust to changes in the quality in initial solutions used;
\item Use of both types of initial solutions would have increased the number of improved solutions. In the table shown in Appendix A the results which are obtained in the two sets of experiments are compared. For 20 benchmark instance with 24 fixed points, the suggested heuristics improved all results in \cite{BassettoM11} obtained with the algorithms. For the results obtained there with a human intervention, 17 results are improved, two results are the same, and for only one ``difficult" instance (instance 52) the obtained result was 0.03\% worse than the result in \cite{BassettoM11}. By using random initial solutions with heuristic $H_{6,4}$ we managed to improve the instance 52 solution (in 37 iterations, 170 seconds), which is now 0.24\% better than the solution in \cite{BassettoM11}. The solution found is shown in Figure~\ref{fig:Instance52}.

\end{itemize}

\section{Summary}  \label{sec:sum} In this paper we have considered the balanced two period TSP. For this  ${\cal NP}$-hard problem, we have described a new polynomially solvable case  which results from the restriction to
Kalmanson matrices. Motivated by this solvable case we suggested  a
simple heuristic. The heuristic managed to improve more than one third of the solutions for the  published instances.

In our approach to the special case and in the simple heuristic we made use
of a  methodological approach ``one sequence - two tours". We formulated then a new heuristic based on the same methodological principle. The core ``engine" in the new heuristic is though different. Instead of DP recursions designed for a special solvable case, we used the DP recursions developed by Held and Karp for the general case of the TSP. To address the computations time and the curse of dimensionality  we suggested an ``assemble-disassemble" heuristic.

Computational experiments on test instances have shown good performance: solutions for 58 out of 60 published benchmark instances have been improved (we refer here to the benchmark results obtained by algorithms in \cite{BassettoM11}).

Our new 2-TSP heuristic has been described by using the terminology for the VRP with two vehicles. We believe that the underlying idea can be used to deal with the case of a larger number of vehicles. 
We have successfully used this approach in real life applications
(to be reported in a forthcoming work) and also successfully took part in 
international contests in logistic optimizations (third prize at the 
2015 VEROLOG competition). 
 We believe that our approach can be used in a variety of VRP settings. This
 will be the subject of future work.

\section*{Acknowledgements}
\nopagebreak  The authors would like to thank referees for valuable comments and suggestions. The authors also thank Tatiana Bassetto and Francesco Mason for providing the benchmark instances.
Financial support  by the University of Warwick DIMAP Center and
 the Austrian Science Fund (FWF): W1230, Doctoral Program ''Discrete Mathematics'' is gratefully acknowledged.

\appendix

\section {Appendix: Results of computational experiments }\nopagebreak[4]

    \begin{tabular}{|c|c|c|c|c|c|c|c|}
    \multicolumn{8}{l}{Results for the instances with {\bf 24} nodes visited twice (48 initial solutions)} \bigstrut[b]\\
    \hline
    \multicolumn{1}{|r|}{\multirow{2}[4]{*}{Instance}} & \multicolumn{2}{c|}{Solutions from [5]} & KS $\times$48 & \multicolumn{2}{c|}{[KS+$H_{5,3}]\times$36} & \multicolumn{2}{c|}{[RP+$H_{5,3}$]} \bigstrut\\
\cline{2-8}        & PC  & PC/h & length & time (s) & length & time (s) & length \bigstrut\\
    \hline
    $I_{41}$ & 30253 & 30147 & 30666 & 71  & 30349 & 91 & \underline{30062} \bigstrut[t]\\
    $I_{42}$ & 33008 & 32020 & 33106 & 93  & \underline{31754} & 92 & 31798 \\
    $I_{43}$ & 31500 & 31500 & 31920 & 96  & 31233 & 90 & \underline{31194} \\
    $I_{44}$ & 30313 & 30170 & 30351 & 86  & 29882 & 91 & \underline{29839} \\
    $I_{45}$ & 27986 & 27857 & 27974 & 64  & 27843 & 90 & \underline{27780} \\
    $I_{46}$ & 30073 & \dotuline{30013} & 30423 & 88 & 30017 & 90 & \underline{30013} \\
    $I_{47}$ & 32106 & 32106 & 32715 & 96  & \underline{31663} & 90 & 31827 \\
    $I_{48}$ & 31004 & 30942 & 31087 & 71  & \underline{30507} & 91 & 30597 \\
    $I_{49}$ & 33663 & 33185 & 34728 & 159  & 33181 & 91 & \underline{33088} \\
    $I_{50}$ & 31266 & 31266 & 31714 & 97  & \underline{31236} & 90 & 31504 \\
    $I_{51}$ & 33722 & 33627 & 33815 & 98  & 33412 & 90 & \underline{33344} \\
    $I_{52}$ & 32353 & \underline{32280} & 32732 & 99  & 32323 & 90 & 32290 \\
    $I_{53}$ & 33287 & 33200 & 34094 & 104  & 33086 & 91 & \underline{32883} \\
    $I_{54}$ & 31973 & 31600 & 32732 & 96  & 31370 & 91 & \underline{31368} \\
    $I_{55}$ & 33837 & \dotuline{33507} & 34709 & 63  & \underline{33507} & 90 & 33560 \\
    $I_{56}$ & 29696 & 29476 & 29995 & 69  & 29496 & 90 & \underline{29227} \\
    $I_{57}$ & 31954 & 31640 & 32526 & 104  &\underline{31427} & 90 & 31432 \\
    $I_{58}$ & 30705 & 30246 & 31257 & 63  & 30487 & 91 & \underline{30183} \\
    $I_{59}$ & 31549 & 31549 & 31650 & 96  & \underline{31471} & 92 & 31591 \\
    $I_{60}$ & 32384 & 32317 & 32794 & 100  & 32193 & 90 & \underline{32137} \bigstrut[b]\\
    \hline
    \end{tabular}%


\section {Appendix: On implementation of the sliding subset approach } 
\begin{figure}[thb]
\unitlength=0.8cm
\begin{center}
{
\begin{picture}(12,9.5)
\put(-2.5,0)
{
\includegraphics[scale=0.8]{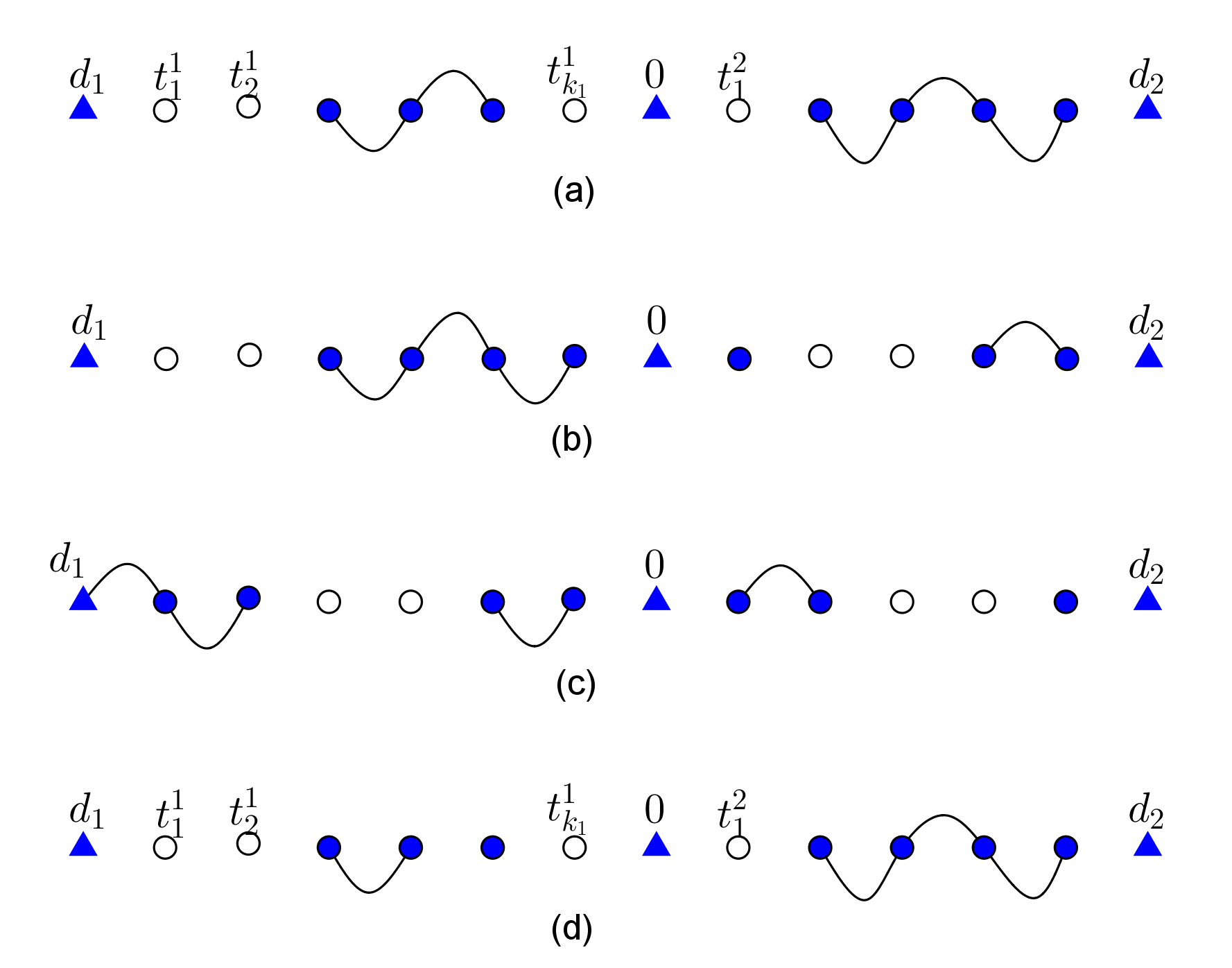}
}
\end{picture} 
}
\end{center}
\caption{ Illustrations of disassembling: (a) first step: $S_1$ and $S_2$ are separated by one sub-path only; (b) an example when $S_1$ and $S_2$ are separated by two sub-paths and a depot; (c) depot and the first sub-path are considered as one customer; (d) modified first step: a sub-path between $S_1$ and $S_2$ (case (a)) is partitioned into a sub-path and a single customer to keep the size of the new 2-VRP fixed.
}
\label{fig:disassembling2}
\end{figure}

\nopagebreak[4]
Figure~\ref{fig:disassembling2} illustrates the various steps 
(and possibilities) of the disassembling process for the case $s=2$ and $l=2$.

 Fig.~\ref{fig:disassembling2}(a) illustrates the first step of the disassembling process as was described above. Sets $S_1$ and $S_2$ are separated by one sub-path in this case. The size of the small 2-VRP  is $2s+5$.

Fig.~\ref{fig:disassembling2}(b) illustrates the outcome of disassembling the route on the next iteration. Notice that depot $0$ is always treated as a separate customer, therefore setting $l=2$ yields the position of $S_2$ as shown in the figure.  There are two sub-paths between sets $S_1$ and $S_2$, and the size of the small 2-VRP  is $2s+6$. 

Consider the step when $S_1$ and $S_2$ are chosen as shown in Fig.~\ref{fig:disassembling2}(c). If the first sub-path did not contain the depot, the size of the problem would have been $2s+7$.  
It means that we can end up with instances with $2s+5$, $2s+6$, and $2s+7$ customers. 

For the implementation, it was convenient to keep the size of the instances fixed at $2s+6$. Therefore we decided (1) to ``glue" the first sub-path with the depot and define it as the depot in the new problem, as shown in Fig.~\ref{fig:disassembling2}(c); 
(2) in case when subsets  $S_1$ and $S_2$ are separated by a single path (for
example, in the first step of disassembling), the last node in the sub-path 
is considered as a sub-path with one node: in this case the instance 
with $2s+5$ customers becomes an instance with $2s+6$ customers (compare Fig.~\ref{fig:disassembling2}(a) and Fig.~\ref{fig:disassembling2}(d)). 

\newpage



\end{document}